\documentclass[12pt]{amsart}
\usepackage{a4wide}
\usepackage[dvips]{epsfig}
\usepackage{amscd}
\usepackage{amssymb}
\usepackage{amsthm}
\usepackage{amsmath}
\usepackage{latexsym}
\usepackage[linktocpage,
                  menucolor=blue,
                  linkcolor=blue,
                  citecolor=red,
                  colorlinks=true]{hyperref}

\theoremstyle{plain}
\newtheorem{theorem}{Theorem}[section]
\theoremstyle{plain}

\newtheorem{corollary}[theorem]{Corollary}

\theoremstyle{definition}

\newtheorem*{remark}{Remark}

{%
\setcounter{enumi}{0}

\begin{enumerate}}%
{\end{enumerate} }
{%
\setcounter{enumi}{0}

\begin{enumerate}}%
{\end{enumerate} }

\numberwithin{equation}{section}
 \allowdisplaybreaks

\title[Optimal portfolios, consumption and  insurance]
 {Optimal investment-consumption and life insurance selection problem under inflation. A BSDE approach }

\date{\today}

\begin{document}

\author{Calisto Guambe}
\address{Department of Mathematics and Applied Mathematics, University of Pretoria, 0002, South Africa}
\address{Department of Mathematics and Informatics, Eduardo Mondlane University, 257, Mozambique}

\email{calistoguambe@yahoo.com.br}

\author{ Rodwell Kufakunesu }

\address{Department of Mathematics and Applied Mathematics
, University of Pretoria, 0002, South Africa}

\email{rodwell.kufakunesu@up.ac.za}


\begin{abstract}
We discuss an optimal investment, consumption and insurance problem of a wage earner under inflation. Assume a wage earner investing in a real money account and three asset prices, namely: a real zero coupon bond, the inflation-linked real money account and a risky share described by jump-diffusion processes. Using the theory of quadratic-exponential backward stochastic differential equation (BSDE) with jumps approach, we derive the optimal strategy for the two typical utilities (exponential and power) and the value function is characterized as a solution of BSDE with jumps. Finally, we derive the explicit solutions for the optimal investment in both cases of exponential and power utility functions for a diffusion case.
\end{abstract}

\maketitle
\section{Introduction}

The problem of asset allocation with life insurance consideration is of great interest to the investor because it protects their dependents if a premature death occurs. Since the optimal portfolio, consumption and life insurance problem by Richard \cite{Richard} in 1975, many works in this direction have been reported in the literature. (See, e.g., Pliska and Ye \cite{Pliska}, Guambe and Kufakunesu \cite{guambe}, Han and Hu \cite{Han}, among others).

In this paper, we discuss an optimal investment, consumption and life insurance problem using the backward stochastic differential equations (BSDE) with jumps approach. Unlike the dynamic programming approach applied in Han and Hu \cite{Han}, this approach allows us to solve the problem in a more general non-Markovian case. For more details on the theory of BSDE with jumps, see e.g., Delong \cite{Delong}, Cohen and Elliott \cite{cohen}, Morlais \cite{morlais2010}, and references therein.
Our results extend, for instance, the paper by  Cheridito and Hu \cite{cheridito} to a jump diffusion setup and we allow the presence of life insurance and inflation risks. Inflation is described as a percentage change of a particular reference index. The inflation-linked products may be used to protect the future cash flow of the wage earner against inflation, which occurs from time to time in some developing economies. Therefore, it make sense to model the inflation-linked products using jump-diffusion processes. For more details on the inflation-linked derivatives, see e.g., Tiong \cite{Tong}, Mataramvura \cite{Mataramvura} and references therein. We consider a model described by a risk-free asset, a real zero coupon bond, an inflation-linked real money account and a risky asset under jump-diffusion processes. These type of processes are more appropriate for modeling the response to some important extreme events that may occur since they allow capturing some sudden changes in the price evolution, as well as, the consumer price index that cannot be explained by models driven by Brownian information. Such events happen due to many reasons, for instance, natural disasters, political situations, etc.

The corresponding quadratic-exponential BSDE with jumps relies on the results by Morlais \cite{morlais2010},  Morlais \cite{morlais2009}, where the existence and uniqueness properties of the quadratic-exponential BSDE with jumps have been proved. Thus, we are also extending the utility maximization problem in Morlais \cite{morlais2009} by including consumption and life insurance. Similar works include Hu {\it et. al.} \cite{hu}, Xing \cite{xing}, Siu \cite{siu}, \O ksendal and Sulem \cite{oksendal2011}, among others.

This paper is organized as follows: in Section 2, we introduce the inflation risks and the related assets: the real zero coupon bond, the inflation-linked real money account, and the risky asset. We also introduce the insurance market and we state the main problem under study. Section 3 is the main section of this paper, we present the general techniques of the BSDE approach and we prove the main results in the exponential and power utility function. Finally, in Section 4, we give some concluding remarks.

\section{Model formulation}
Suppose we have a wage earner investing in a finite investment period $T<\infty$, which can be interpreted as a retirement time. Consider a complete filtered probability space ($\Omega,\mathcal{F},\{\mathcal{F}_t\}_{0\leq t\leq T},\mathbb{P}$), where
$\{\mathcal{F}_t\}_{t\in[0,T]}$ is a filtration satisfying the usual conditions (Protter \cite{Protter}). Denote by $W_r$ and $W_n$ the Brownian motions underlying the risks driven by the real and nominal term structures. We also define the Brownian motions $W_I$ and $W_S$, the drivers in the inflation rate and the risky asset. We assume that $W_r,\,W_n,\,W_I$ and $W_S$ are independent. Note that if we allow the correlations among $W_r,\,W_n,\,W_I$ and $W_S$, i.e., $dW_k(t)dW_I(t)=\rho_{kI}dt$; $dW_k(t)dW_S(t)=\rho_{kS}dt$ for $k\in\{r,n\}$ and $dW_I(t)dW_S(t)=\rho_{IS}dt$, where $\rho_{ij}$ are the correlation coefficients, may result in a highly nonlinear BSDE with jumps which the existence and uniqueness of its solution has not yet been established. Moreover, we consider a Poisson process $N$ independent of  $W_r,\, W_n,\,W_I$ and $W_S$, associated with the complete filtered probability space $(\Omega,\mathcal{F},\{\mathcal{F}_t\},\mathbb{P})$ with the intensity measure $dt\times d\nu(z)$, where $\nu$ is the $\sigma$-finite Borel measure on $\mathbb{R}\setminus\{0\}$. A $\mathbb{P}$-martingale compensated Poisson random measure is given by:
\begin{equation*}
    \tilde{N}(dt,dz):=N(dt,dz)-\nu(dz)dt\,.
\end{equation*}
Furthermore, we consider the following spaces:
\begin{itemize}
  \item $\mathbb{L}^2(\mathbb{R})$- the space of random variables
$\xi:\Omega\mapsto\mathbb{R}$, such that $\mathbb{E}[\,|\xi|^2]<\infty$.
  \item $\mathbb{H}^2(\mathbb{R})$- the space of measurable functions
$Z:\mathbb{R}\mapsto\mathbb{R}$ such that
$$\mathbb{E}\left[\int_{\mathbb{R}}|Z(t)|^2dt\right]<\infty\,.$$
  \item $\mathbb{S}^2(\mathbb{R})$- the space of adapted c\`adl\`ag processes
$Y:\Omega\times[0,T]\mapsto\mathbb{R}$ such that
$$\mathbb{E}[\sup|Y(t)|^2]<\infty$$ and
  \item $\mathbb{H}_\nu^2$- the space of predictable processes
$\Upsilon:\Omega\times[0,T]\times\mathbb{R}\mapsto\mathbb{R}$, such that
\begin{equation*}
\mathbb{E}\left[\int_0^T\int_{\mathbb{R}}|\Upsilon(t,z)|^2\nu(dz)dt\right]
<\infty.
\end{equation*}

\end{itemize}

Let $r$ denote the real and $n$ the nominal forward rates, defined for $k\in\{r,n\}$, by:
\begin{equation*}
    f_k(t,T):=f_k(0,T)+\int_0^t\alpha_k(s,T)ds+\int_0^t\sigma_k(s,T)dW_k(s)+\int_0^t\gamma_k(s,T,z)\tilde{N}(ds,dz)\,,
\end{equation*}
where the coefficients $\alpha_k(t,T),\ \ \sigma_k(t,T)$ and $\gamma_k(t,T,z)$ are $\mathcal{F}_t$-predictable bounded processes, satisfying the following condition:
$$
\int_0^T\left[|\alpha_k(t,T)|+\sigma_k^2(t,T)+\int_{\mathbb{R}}\gamma_k^2(t,T,z)\nu(dz)\right]dt<\infty, \ \ \ \rm{a.s.}
$$
We denote by $r_k(t)=f_k(t,t)$ the corresponding spot rate at time $t$.

It is well known that the price of the real (nominal) bond is given by

$$
P_k(t,T):=\exp\left\{-\int_t^Tf_k(t,s)ds\right\}\,.
$$
An application of the It\^o's formula yields:

\begin{equation*}
    dP_k(t,T)=P_k(t,T)\left\{a_k(t,T)dt+b_k(t,T)dW_k(t)+\int_{\mathbb{R}}c_k(t,T,z)\tilde{N}(dt,dz)\right\}\,,
\end{equation*}
where
$$
b_k(t,T):=-\int_t^T\sigma_k(t,s)ds\,; \ \ \ c_k(t,T,z):=-\int_t^T\gamma_k(t,s,z)ds
$$
and
$$
   a_k(t,T):=r_k(t)-\int_t^T\alpha_k(t,s)ds+\frac{1}{2}\|b_k(t,T)\|^2-\int_{\mathbb{R}}c_k(t,T,z)\nu(dz)\,.
$$

We suppose the existence of an inflation index $I(t)$, i.e., the consumer price index (CPI) governed by the following stochastic differential equation (SDE)
\begin{equation*}
    dI(t)=I(t)\left[\mu_I(t)dt+\sigma_I(t)dW_I(t)+\int_{\mathbb{R}}\gamma_I(t,z)\tilde{N}(dt,dz)\right]\,,
\end{equation*}
where the expected inflation rate $\mu_I(t)$, the volatility $\sigma_I(t)$ and the dispersion rate $\gamma_I(t,z)>-1$ are $\mathcal{F}_t$-predictable bounded processes, satisfying the following integrability condition
$$
\int_0^T\left[|\mu_I(t)|+\sigma_I^2(t)+\int_{\mathbb{R}}\gamma_I^2(t,z)\nu(dz)\right]dt<\infty, \ \ \ \rm{a.s.}
$$

The financial market consists of four assets, namely a real (nominal) money account $B_k(t)$ defined by
$$
B_k(t)=\exp\left\{\int_0^tr_k(s)ds\right\}\,.
$$
A real zero coupon bond price $P_r^*(t,T)$ defined as
$$
P_r^*(t,T)=I(t)P_r(t,T).
$$
Applying the It\^o's product rule, we have that
\begin{eqnarray*}
  dP_r^*(t,T) &=& I(t)dP_r(t,T)+P_r(t,T)dI(t)+d[I(t),P_r(t,T)] \\
   &=& P_r^*(t,T)\left[\tilde{A}(t,T)dt+b_r(t,T)dW_r(t)+\sigma_I(t)dW_I(t)+\int_{\mathbb{R}}\tilde{C}(t,T,z)\tilde{N}(dt,dz)\right],
\end{eqnarray*}
where
$$
\tilde{A}(t,T):=a_r(t,T)+\mu_I(t)+\int_{\mathbb{R}}c_r(t,T,z)\gamma_I(t,z)\nu(dz)
$$
and
$$
\tilde{C}(t,T,z):=c_r(t,T,z)+\gamma_I(t,z)+c_r(t,T,z)\gamma_I(t,z)\,.
$$

We also define the inflation-linked real money account $B_r^*(t)$ by
$$
B_r^*(t):=I(t)B_r(t).
$$
Then, by It\^o's formula, we can easily see that it is governed by the following SDE:
$$
dB_r^*(t)=B_r^*(t)\left[(r_r(t)+\mu_I(t))dt+\sigma_I(t)dW_I(t)+\int_{\mathbb{R}}\gamma_I(t,z)\tilde{N}(dt,dz)\right]\,.
$$

Finally, we define the risky asset by the following geometric jump-diffusion process
$$
dS(t)=S(t)\left[\mu_S(t)dt+\sigma_S(t)dW_S(t)+\int_{\mathbb{R}}\gamma_S(t,z)\tilde{N}(dt,dz)\right]\,,
$$
where the mean rate of return $\mu_S(t)$, the volatility $\sigma_S(t)$ and the dispersion rate $\gamma_S(t,z)>-1$ are $\mathcal{F}_t$-predictable bounded processes, satisfying the following integrability condition
$$
\int_0^T\left[|\mu_S(t)|+\sigma_S^2(t)+\int_{\mathbb{R}}\gamma_S^2(t,z)\nu(dz)\right]dt<\infty, \ \ \ \rm{a.s.}
$$

For later use, we define the following processes (also called market price of risks) $\varphi_1:=\frac{\tilde{A}-r_r}{b_r}$, $\varphi_2:=\frac{\mu_I}{\sigma_I}$ and $\varphi_3:=\frac{\mu_S-r_r}{\sigma_S}$,  provided that $b_r, \sigma_I,\sigma_S\neq0$.

As in Guambe and Kufakunesu \cite{guambe}, we consider a wage earner whose lifetime is a
nonnegative random variable $\tau$ defined on the probability space $(\Omega,
\mathcal{F},\mathbb{P})$.
Consider that $\tau$ has a probability density function $g(t)$ and the distribution function is given by
$$
G(t):=\mathbb{P}(\tau<t)=\int_0^tg(s)ds\,.
$$
The probability that the life time $\tau>t$ is given by
$$
\bar{G}(t):=\mathbb{P}(\tau\geq t\mid\mathcal{F}_t)=1-G(t)\,.
$$
The instantaneous force of mortality $\lambda(t)$ for the policyholder to be alive at time $t$ is defined by
\begin{eqnarray*}
  \lambda(t) &:=& \lim_{\Delta t\rightarrow0}\frac{\mathbb{P}(t\leq\tau<t+\Delta t|\tau\geq t)}{\Delta t} \\
   &=& \lim_{\Delta t\rightarrow0}\frac{\mathbb{P}(t\leq\tau<t+\Delta t)}{\Delta t\mathbb{P}(\tau\geq t)} \\
   &=&  \frac{1}{\bar{G}(t)}\lim_{\Delta t\rightarrow0}\frac{G(t+\Delta t)-G(t)}{\Delta t} \\
   &=& \frac{g(t)}{\bar{G}(t)}=-\frac{d}{dt}(\ln(\bar{G}(t)))\,.
\end{eqnarray*}

 Then, the conditional survival probability of the policyholder is given by
 \begin{equation}\label{survival}
\bar{G}(t)=\mathbb{P}(\tau>t|\mathcal{F}_t)=\exp\left(-\int_0^t\lambda(s)ds\right),
 \end{equation}
 and the conditional survival probability density of the death of the policyholder by
 \begin{equation}\label{death}
    g(t):=\lambda(t)\exp\left(-\int_0^t\lambda(s)ds\right).
 \end{equation}

We suppose existence of an insurance market, where the term life insurance is continuously traded. We assume that the wage earner is paying premiums at the rate $p(t)$, at time $t$ for the life insurance contract and the insurance company will pay $p/\eta(t)$ to the beneficiary for his death, where the $\mathcal{F}_t$-adapted process $\eta(t)>0$ is the premium insurance ratio.
When the wage earner dies, the total legacy to his beneficiary is given by
 \begin{equation*}
    \ell(t):=X(t)+\frac{p(t)}{\eta(t)},
 \end{equation*}
where $X(t)$ is the wealth process of the wage earner at time $t$ and $p(t)/\eta(t)$ the insurance benefit paid by the insurance company to the beneficiary if death occurs at time $t$.

Let $c(t)$ be the consumption rate of the wage earner and $\theta(t):=(\theta_1(t),\theta_2(t),\theta_3(t))$ be the vector of the amounts of the wage earner's wealth invested in the real zero coupon bond $P_r^*$, the inflation-linked real money account $B_r^*$ and the risky asset $S$ respectively, satisfying the following integrability condition.
\begin{equation}\label{integrabilitycp}
\int_0^T\Bigl[c(t)+p(t)+\sum_{i=1}^3\theta_i^2(t)\Bigl]dt<\infty, \ \ \ \rm{a.s.}
\end{equation}

Furthermore, we assume that the shares are divisible, continuously traded and there are no transaction costs, taxes or short-selling constraints in the trading. Then the wealth process $X(t)$ is defined by the following (SDE):

\begin{eqnarray}\label{wealth}
  dX(t) &=& [r_r(t)X(t)+\langle\theta(t),\hat{\mu}(t)\rangle-c(t)-p(t)]dt+\theta_1(t)b_r(t,T)dW_r(t) \\ \nonumber
   && +(\theta_1(t)+\theta_2(t))\sigma_I(t)dW_I(t) +\theta_3(t)\sigma_S(t)dW_S(t) \\ \nonumber
   &&  +\int_{\mathbb{R}}\langle\theta(t),\hat{\gamma}(t,T,z)\rangle\tilde{N}(dt,dz)\,, \ \ \ t\in[0,\tau\wedge T]
\end{eqnarray}
where $\hat{\mu}(t):=(\tilde{A}(t,T)-r_r(t),\mu_I(t),\mu_S(t)-r_r(t))$, $\hat{\gamma}(t,T,z):=(\tilde{C}(t,T,z),\gamma_I(t,z),\gamma_S(t,z))$, $\tau\wedge T:=\min\{\tau;T\}$ and $\langle\cdot,\cdot\rangle$ is the inner product in $\mathbb{R}^n$.

The wage earner faces the problem of choosing the optimal strategy $\mathcal{A}:=\{(\theta,c,p):=(\theta(t),c(t),p(t))_{t\in[0,T]}$\} which maximizes the discounted expected utilities from the consumption during his/her lifetime $[0,\tau\wedge T]$, from the wealth if he/she is alive until the terminal time $T$ and from the legacy if he/she dies before time $T$. Suppose that the discount process rate $\varrho(t)$ is positive and
$\mathcal{F}_t$-adapted process. This problem can be defined by the following performance functional (for more details see, e.g., Pliska and Ye \cite{Pliska}, Oksendal and Sulem \cite{Oksendal}, Guambe and Kufakunesu \cite{guambe}).

\begin{eqnarray}\label{maximumutility}
   J(0,x_0;\theta,c,p) &:=&
ess\sup_{(\theta,c,p)\in\mathcal{A}}\mathbb{E}\Bigl[\int_0^{\tau\wedge T}
e^{-\int_0^s\varrho(u)du}U(c(s))ds  \\ \nonumber &&
+e^{-\int_0^{\tau}\varrho(u)du}U(\ell(\tau))\mathbf{1}_{\{\tau\leq T\}}
+e^{-\int_0^T\varrho(u)du}U(X(T))\mathbf{1}_{\{\tau>T\}}\Bigl],
\end{eqnarray}
where $\mathbf{1}_A$ is a characteristic function defined on a set $A$ and $U$ is the utility function for the consumption, legacy and terminal wealth.

Note that from the conditional survival probability of the wage earner \eqref{survival} and the conditional survival probability density of death of the wage earner \eqref{death}, we can write a dynamic version of the functional \eqref{maximumutility} by:
\begin{eqnarray}\nonumber
  J(t) &=&  \mathbb{E}_{t,x}\Bigl[\int_t^T
e^{-\int_t^s(\varrho(u)+\lambda(u))du}[U(c(s)) +\lambda(s)U(\mathcal{\ell}(s))]ds \\ \label{functional0}
   && \ \ \ \ \ \  +e^{-\int_t^T(\varrho(u)+\lambda(u))du}U(X(T))\mid\mathcal{F}_t\Bigl].
\end{eqnarray}
Thus, the problem of the wage earner is to maximize the above dynamic performance functional under the admissible strategy $\mathcal{A}$. Therefore, the value function  $V(t,x)$ can be restated in the following form:
\begin{equation}\label{valuefunction}
V(t,x)=ess\sup_{(\theta,c,p)\in\mathcal{A}}J(t,x,\theta,c,p)\,.
\end{equation}

The set of strategies $\mathcal{A}:=\{(\theta,c,p):=(\theta(t),c(t),p(t))_{t\in[0,T]}$\} is said to be admissible if the SDE \eqref{wealth} has a unique strong solution such that $X(t)\geq0$, $\mathbb{P}$-a.s. and

\begin{eqnarray*}
   &&
\mathbb{E}_{t,x}\Bigl[\int_t^T
e^{-\int_t^s(\varrho(u)+\lambda(u))du}[U(c(s)) +\lambda(s)U(\mathcal{\ell}(s))]ds \\
   && \ \ \ \ \ \  +e^{-\int_t^T(\varrho(u)+\lambda(u))du}U(X(T))\mid\mathcal{F}_t\Bigl]<\infty\,.
\end{eqnarray*}
The nonnegative condition of the wealth process contains a non-borrowing constraints that prevents the family from borrowing for consumption and life insurance at any time $t\in[0,T]$.

In order to solve our optimization using the quadratic-exponential BSDE's with jumps and to make the proofs easier, we introduce, in addition to the integrability condition \eqref{integrabilitycp}, the constraints in the admissible strategy $\mathcal{A}$ as follows: let $\mathcal{C}\subset\mathcal{P}$ and $\mathcal{D}\subset\mathcal{P}$, where $\mathcal{P}$ denotes the set of real valued predictable processes $(c(t))_{0\leq t\leq T}$, $(p(t))_{0\leq t\leq T}$, and $\mathcal{Q}\subset\mathcal{P}^{1\times 3}$, where $\mathcal{P}^{1\times 3}$ represents the set of all predictable processes $(\theta_1(t),\theta_2(t),\theta_3(t))_{0\leq t\leq T}$. In the exponential case, we assume that the admissible strategy $(c(t),p(t),\theta(t))\in\mathcal{C}\times\mathcal{D}\times\mathcal{Q}$. For the power utility case, the consumption, investment and life insurance strategies will be denoted by their fractions of the total wealth, that is, $c=\xi X$, $\theta=\pi X$, and $p=\zeta X$. We assume that $(\xi(t),\zeta(t),\pi(t))\in\mathcal{C}\times\mathcal{D}\times\mathcal{Q}$.

We assume that $\mathcal{C},\,\mathcal{D}$ and $\mathcal{Q}$ are closed and compact sets.

Moreover, we introduce the notion of martingales of bounded mean oscillation ({\it BMO}-martingales) for jump-diffusion processes as in Morlais \cite{morlais2009}. We say that a Martingale $M$ is in the class of {\it BMO}-martingales if there exists a constant $K>0$, such that, for all $\mathcal{F}$-stopping times $\mathcal{T}$,
$$
ess\sup_{\Omega}\mathbb{E}[[M]_T-[M]_{\mathcal{T}}\mid\mathcal{F}_{\mathcal{T}}]\leq K^2 \ \ \ \ \ \rm{and} \ \ \ \ |\Delta M_{\mathcal{T}}|\leq K^2\,,
$$
where $[M]$ denotes a quadratic variation of a process $M$. For the diffusion case, the $BMO$-martingale property follows from the first condition, whilst in a jump-diffusion case, we need to ensure the boundedness of the jumps of the local martingale $M$.

\section{The BSDE approach to optimal investment, consumption and insurance}

In this section, we solve the optimal investment, consumption and life insurance problem under inflation using the BSDE with jumps approach. For more details on the theory of BSDEs with jumps see, e.g., Delong \cite{Delong}. We then consider two utility functions, namely, the exponential utility and the power utility. The techniques we use  are similar to Morlais \cite{morlais2009}, Cheridito and Hu \cite{cheridito}, Xing \cite{xing}.

Define the following BSDE with jumps:
\begin{eqnarray}\label{bsde}
  dY(t) &=& -h(t, Y(t), Z_1(t), Z_2(t), Z_3(t), \Upsilon(t,\cdot), \theta(t), c(t),p(t))dt+Z_1(t)dW_r(t) \\ \nonumber
   && +Z_2(t)dW_I(t)+Z_3(t)dW_S(t)+\int_{\mathbb{R}}\Upsilon(t,z)\tilde{N}(dt,dz); \\ \nonumber
  Y(T) &=& 0\,.
\end{eqnarray}

The aforementioned approach is based on the following:
Consider the process
\begin{eqnarray}\nonumber
  \mathcal{R}(t) &=&  \int_0^t
e^{-\int_0^s(\varrho(u)+\lambda(u))du}[U(c(s)) +\lambda(s)U(\mathcal{\ell}(s))]ds \\ \nonumber
   && \ \ \ \ \ \  +e^{-\int_t^T(\varrho(u)+\lambda(u))du}U(X(t)-Y(t))\,,
\end{eqnarray}
with the initial condition $\mathcal{R}(0)=U(x-Y(0))$. Here, $X(t)$ represents the wealth process \eqref{wealth} and $Y(t)$ part of the solution $(Y,Z_1,Z_2,Z_3,\Upsilon)$ of the BSDE with jumps \eqref{bsde}.
Applying the generalized It\^o's formula, we have

\begin{eqnarray*}
 d\mathcal{R}(t)  &=&  e^{-\int_0^t(\varrho(s)+\lambda(s))ds}\Bigl\{[\Lambda(t,y,z_1,z_2, z_3,\upsilon,\theta,c,p) +h(t,y,z_1,z_2, z_3,\upsilon)]dt \\
   && U'(X(t)-Y(t))[(\theta_1(t)b_r(t,T)+z_1)dW_r(t) \\
&& +((\theta_1(t)+\theta_2(t))\sigma_I(t)+z_2)dW_I(t)+(\theta_3(t)\sigma_S(t)+z_3)dW_S(t)] \\
&& +\int_{\mathbb{R}}[U(X(t)-Y(t)+\langle\theta(t),\hat{\gamma}(t,T,z)\rangle+\upsilon(t,z))-U(X(t)-Y(t))]\tilde{N}(dt,dz)\Bigl\}\,,
\end{eqnarray*}
where
\begin{eqnarray*}
 && \Lambda(t,y,z_1,z_2, z_3,\upsilon,\theta,c,p) \\
&=& -\Bigl\{[U(c(t))+\lambda(t)U(\mathcal{\ell}(t)) -(\varrho(t)+\lambda(t))U(X(t)-Y(t))] \\
   && +U'(X(t)-Y(t))[r_r(t)X(t)+\langle\theta(t),\hat{\mu}(t)\rangle-c(t)-p(t)] \\
   &&  + \frac{1}{2}U''(X(t)-Y(t))[(\theta_1(t)b_r(t,T)+z_1)^2+ ((\theta_1(t)+\theta_2(t))\sigma_I(t)+z_2)^2 \\
&& + (\theta_3(t)\sigma_S(t)+z_3)^2 ] +\int_{\mathbb{R}}[U(X(t)-Y(t)+\langle\theta(t),\hat{\gamma}(t,T,z)\rangle+\upsilon(t,z)) \\
&& -U(X(t)-Y(t)) -U'(X(t)-Y(t))(\langle\theta(t),\hat{\gamma}(t,T,z)\rangle+\upsilon(t,z))]\nu(dz)\Bigl\}\,.
\end{eqnarray*}
Note that we can write $\mathcal{R}$ as

\begin{eqnarray}\nonumber
    \mathcal{\mathcal{R}}(t) &=& \mathcal{R}(0)+\int_0^te^{-\int_0^s(\varrho(u)+\lambda(u))du}\Bigl\{[\Lambda(s,y,z_1,z_2, z_3,\upsilon,\theta,c,p) +h(s,y,z_1,z_2, z_3,\upsilon)]ds \\ \label{martingale}
&& + \{a \ \ local \ \ martingale\}\,.
\end{eqnarray}

We define the generator $h$ by
$$
h(s,y,z_1,z_2, z_3,\upsilon)=\inf_{(\theta,c,p)}\Lambda(s,y,z_1,z_2, z_3,\upsilon,\theta,c,p)\,.
$$
Then we can see that \eqref{martingale} is decreasing, hence $\mathcal{R}$ is a local super-martingale and we can choose a strategy $(\theta^*,c^*,p^*)$ such that the drift process in \eqref{martingale} is equal to zero, therefore, $\mathcal{R}$ is a local martingale and prove that $(\theta^*,c^*,p^*)$ is the optimal strategy.

We will establish the existence and uniqueness properties of the solution $(Y,Z_1,Z_2,Z_3,\Upsilon) \in\mathbb{S}(\mathbb{R})\times\mathbb{H}^2(\mathbb{R})\times\mathbb{H}^2(\mathbb{R}) \times\mathbb{H}^2(\mathbb{R})\times\mathbb{H}^2_\nu(\mathbb{R})$ of the BSDE with jumps \eqref{bsde}, as well as the characterization of the optimal strategy $(\theta^*,c^*,p^*)$, for the specific utilities in the following subsections.
~~\\

\subsection{The exponential utility}
~~\\

We consider the exponential utility function of the form
\begin{equation}\label{exponential}
U(x)=-e^{-\delta x}, \ \ \delta>0\,.
\end{equation}

The functional \eqref{functional0}, is then given by
\begin{eqnarray}\nonumber
  J(t) &=&  -\mathbb{E}_{t,x}\Bigl[\int_t^T
e^{-\int_t^s(\varrho(u)+\lambda(u))du}[e^{-\delta c(s)} +\lambda(s)e^{-\delta \ell(s)}]ds \\ \label{functionalexp}
   && \ \ \ \ \ \  +e^{-\int_t^T(\varrho(u)+\lambda(u))du}\cdot e^{-\delta X_t(T)}\mid\mathcal{F}_t\Bigl]\,.
\end{eqnarray}

We then state the main result of this subsection
\begin{theorem}\label{theorem-power}
Suppose that the utility function is given by \eqref{exponential}. Then the optimal value function of the optimization problem \eqref{valuefunction} is given by
\begin{equation}\label{optimalwealth}
V(t,x)=-\exp(-\delta(x-Y(t))),
\end{equation}
where $Y(t)$ is part of the solution $(Y,Z_1,Z_2,Z_3,\Upsilon)$ of the BSDE with jumps \eqref{bsde}, with terminal condition $Y(T)=0$ and the generator $h$ given by

\begin{eqnarray}\label{generatorh}
 && h(t, y, z_1, z_2, z_3, \upsilon(\cdot)) \\ \nonumber
&=& (1-r_r(t))X(t)+\frac{1}{\delta}\Bigl(1+\eta(t)-\varrho(t)-\lambda(t)+\ln\delta +\eta(t)\ln\Bigl(\frac{\delta\lambda(t)}{\eta(t)}\Bigl)\Bigl) - \Bigl(1+\frac{\eta(t)}{\delta}\Bigl)y \\ \nonumber
&& +\inf_{\theta}\Bigl\{\frac{\delta}{2}\Bigl[\left|\theta_1(t)b_r(t,T)-\left(z_1+\frac{\varphi_1(t)}{\delta}\right)\right|^2 +\left|(\theta_1(t)+\theta_2(t))\sigma_I(t)-\left(z_2+\frac{\varphi_2(t)}{\delta}\right)\right|^2 \\ \nonumber
   &&   + \left|\theta_3(t)\sigma_S(t)-\left(z_3+\frac{\varphi_3(t)}{\delta}\right)\right|^2\Bigl] \\ \nonumber
&& + \frac{1}{\delta}\int_{\mathbb{R}}[\exp(\delta(\upsilon(t,z)-\langle\theta(t),\hat{\gamma}(t,T,z)\rangle))-1- \delta(\upsilon(t,z)-\langle\theta(t),\hat{\gamma}(t,T,z)\rangle)]\nu(dz)\Bigl\} \\ \nonumber
&& -(\varphi_1(t)z_1+\varphi_2(t)z_2+\varphi_3(t)z_3) -\frac{1}{2\delta}\left(\varphi_1^2(t) +\varphi_2^2(t)+\varphi_3^2(t)\right) \,.
\end{eqnarray}
Furthermore, the optimal admissible strategy $(\theta^*(t),c^*(t),p^*(t))$ is given by
$$
c^*(t)=X^{(\theta^*,c^*,p^*)}(t)-Y(t)+\frac{1}{\delta}\ln\delta; \ \ \ \ p^*(t)=\eta(t)\Bigl[\frac{1}{\delta}\ln\Bigl(\frac{\delta\lambda(t)}{\eta(t)}\Bigl)-Y(t)\Bigl]
$$
and
\begin{eqnarray}\label{optimalstrategy}
 && \theta^*(t) \\ \nonumber
&=& \inf_{\theta}\Bigl\{\frac{\delta}{2}\Bigl[\left|\theta_1(t)b_r(t,T)-\left(z_1+\frac{\varphi_1(t)}{\delta}\right)\right|^2 \\ \nonumber
   &&  +\left|(\theta_1(t)+\theta_2)\sigma_I(t)-\left(z_2+\frac{\varphi_2(t)}{\delta}\right)\right|^2 + \left|\theta_3(t)\sigma_S(t)-\left(z_3+\frac{\varphi_3(t)}{\delta}\right)\right|^2\Bigl] \\ \nonumber
&&\,+ \frac{1}{\delta}\int_{\mathbb{R}}[\exp(\delta(\upsilon(t,z)-\langle\theta(t),\hat{\gamma}(t,T,z)\rangle))-1- \delta(\upsilon(t,z)-\langle\theta(t),\hat{\gamma}(t,T,z)\rangle)]\nu(dz)\Bigl\}\,.
\end{eqnarray}
\end{theorem}

Note that for the optimal investment strategy $\theta^*(t)=(\theta^*_1(t),\theta^*_2(t),\theta^*_3(t))$, the solution \eqref{optimalstrategy} is not explicit. We then obtain an explicit solution for a special case where there is no jumps, that is $\nu=0$. Applying the first order condition of optimality in \eqref{optimalstrategy}, we prove that the optimal strategy $\theta^*(t)=(\theta^*_1(t),\theta^*_2(t),\theta^*_3(t))$ is given by the following corollary.

\begin{corollary}\label{corr-power}
Assume that $\nu=0$, then the optimal portfolio strategy $(\theta^*_1(t),\theta^*_2(t),\theta^*_3(t))$, for all $t\in[0,T]$ is given by
\begin{eqnarray*}\label{corollary-powerp}
  \theta_1^*(t) &=& \frac{\tilde{A}(t,T)-r_r(t)-\mu_I(t)}{\delta  b_r^2(t,T)}+\frac{Z_1(t)}{b_r(t,T)} \\
  \theta_2^*(t) &=& \frac{1}{\delta}\left[\left(\frac{1}{\sigma_I^2(t)}+\frac{1}{b_r^2(t,T)}\right)\mu_I(t) -\frac{\tilde{A}(t,T)-r_r(t)}{b_r^2(t,T)} +\frac{Z_2(t)}{\sigma_I(t)} -\frac{Z_1(t)}{b_r(t,T)}\right] \\
  \theta_3^*(t) &=& \frac{\mu_S(t)-r_r(t)}{\delta\sigma_S^2(t)}+\frac{Z_3(t)}{\sigma_S(t)}\,,
\end{eqnarray*}
where $(Z_1(t), Z_2(t), Z_3(t))$ is part of the solution $(Y, Z_1, Z_2, Z_3)$ of the following BSDE.
\begin{eqnarray*}
  dY(t) &=& -h(t, Y(t), Z_1(t), Z_2(t), Z_3(t), \theta^*(t), c^*(t),p^*(t))dt+Z_1(t)dW_r(t) \\ \nonumber
   && +Z_2(t)dW_I(t)+Z_3(t)dW_S(t); \\ \nonumber
  Y(T) &=& 0\,.
\end{eqnarray*}
\end{corollary}

Before we prove the main theorem of this subsection, we establish the assumptions for the existence and uniqueness solution of a BSDE with quadratic growth. Suppose we are given a BSDE  \eqref{bsde}, with terminal condition $Y(T)=0$ and a generator $h$ given by \eqref{generatorh}. From the boundedness of the associated parameters, there exists a constant $K>0$ such that
\begin{eqnarray}\nonumber
|h(t, y,z_1,z_2,z_3,\upsilon)| &\leq& K\Bigl(1+|y|+|z_1|^2+|z_2|^2+|z_3|^2 \\ \label{boundedness}
&&  +\frac{1}{\delta}\int_{\mathbb{R}}[\exp(\delta(\upsilon(t,z)-\langle\theta(t),\hat{\gamma}\rangle))-1- \delta(\upsilon(t,z)-\langle\theta(t),\hat{\gamma}\rangle)]\nu(dz)\Bigl)\,.
\end{eqnarray}
Moreover,
\begin{equation}\label{growthcondition}
|h(t, y,z_1,z_2,z_3,\upsilon)-h(t, y',z_1',z_2',z_3',\upsilon)|\leq K\Bigl(|y-y'|+\sum_{i=1}^3(1+|z_i|+|z_i'|)|z_i-z_i'|\Bigl)
\end{equation}
and
$$
|h(t, y,z_1,z_2,z_3,\upsilon)-h(t, y,z_1,z_2,z_3,\upsilon')|\leq \int_{\mathbb{R}}\Phi(\upsilon,\upsilon')(\upsilon-\upsilon')\nu(dz)\,,
$$
where
\begin{eqnarray*}
\Phi(\upsilon,\upsilon') &=& \sup_{\theta}\left(\int_0^1m'(s(\upsilon-\langle\theta,\hat{\gamma}\rangle) +(1-s)(\upsilon'-\langle\theta,\hat{\gamma}\rangle)(z))ds\right)\mathbf{1}_{\upsilon\geq\upsilon'} \\
&& + \inf_{\theta}\left(\int_0^1m'(s(\upsilon-\langle\theta,\hat{\gamma}\rangle) +(1-s)(\upsilon'-\langle\theta,\hat{\gamma}\rangle)(z))ds\right)\mathbf{1}_{\upsilon<\upsilon'}\,,
\end{eqnarray*}
for the function $m$ defined by $m(x)=\frac{\exp(\delta x)-1-\delta x}{\delta}$.

Then, it follows from Morlais \cite{morlais2009}, \cite{morlais2010}, that the BSDE with jumps \eqref{bsde}, with terminal condition $Y(T)=0$ and a generator \eqref{generatorh} has a unique solution $(Y,Z_1,Z_2,Z_3,\Upsilon)\in\mathbb{S}(\mathbb{R})\times\mathbb{H}^2(\mathbb{R})\times\mathbb{H}^2(\mathbb{R}) \times\mathbb{H}^2(\mathbb{R})\times\mathbb{H}^2_\nu(\mathbb{R})$.\\

\noindent {\it Proof of Theorem \ref{theorem-power}.}

Define a family of processes

\begin{eqnarray}\label{processr}
  \mathcal{R}_1^{(\theta,c,p)}(t) &=&  -\int_0^t
e^{-\int_0^s(\varrho(u)+\lambda(u))du}[e^{-\delta c(s)} +\lambda(s)e^{-\delta \ell(s)}]ds \\ \nonumber
   && \ \ \ \ \ \  -e^{-\int_0^t(\varrho(u)+\lambda(u))du}\cdot e^{-\delta (X^{(\theta,c,p)}(t)-Y(t))}\,.
\end{eqnarray}
We aim to construct the process $\mathcal{R}_1^{(\theta,c,p)}$ such that for each strategy $(\theta,c,p)\in\mathcal{A}$, it is a super-martingale  and there exists a strategy $(\theta^*,c^*,p^*)\in\mathcal{A}$ such that $\mathcal{R}_1^{(\theta^*,c^*,p^*)}$ is a martingale.

Applying the It\^o's formula for the process \eqref{processr}, we have

\begin{eqnarray}\label{formularexp}
 && d\mathcal{R}_1(t) \\ \nonumber
&=& \delta e^{-\int_0^t(\varrho(u)+\lambda(u))du}\cdot e^{-\delta (X(t)-Y(t))} \Bigl\{\Bigl[ -e^{\delta (X(t)-Y(t))}\left[e^{-\delta c(t)} +\lambda(t)e^{-\delta \ell(t)}\right]  \\ \nonumber
&& +\frac{1}{\delta}(\varrho(t)+\lambda(t))+ r_r(t)X(t)+\langle\theta(t),\hat{\mu}(t)\rangle-c(t)-p(t) +h(t,z_1,z_2,z_3,\upsilon(\cdot)) \\ \nonumber
   && -\frac{\delta^2}{2} [(z_1-\theta_1b_r(t,T))^2+(z_2-(\theta_1(t)+\theta_2(t))\sigma_I(t))^2 +(z_3-\theta_3(t)\sigma_S(t))^2] \\ \nonumber
&& + \frac{1}{\delta}\int_{\mathbb{R}}[1-\delta(\langle\theta(t),\hat{\gamma}(t,T,z)\rangle-\upsilon(t,z)) -\exp(-\delta(\langle\theta(t),\hat{\gamma}(t,T,z)\rangle-\upsilon(t,z)))]\nu(dz)\Bigl]dt \\ \nonumber
&& -(z_1-\theta_1(t)b_r(t,T))dW_r(t)-(z_2-(\theta_1(t)+\theta_2(t)))\sigma_I(t)dW_I(t)  \\ \nonumber
&& -(z_3-\theta_3(t)\sigma_S(t))dW_S(t) + \int_{\mathbb{R}}[1-\exp(-\delta(\langle\theta(t),\hat{\gamma}(t,T,z)\rangle-\upsilon(t,z)))]\tilde{N}(dt,dz)\Bigl\}
\end{eqnarray}
with the initial condition $\mathcal{R}_1(0)=-\exp(-\delta(x-Y(0)))$.

Note that the drift process of the family $\mathcal{R}_1$ is given by
\begin{eqnarray}
 && A(t) \\ \nonumber
&=& \delta e^{-\int_0^t(\varrho(u)+\lambda(u))du}\cdot e^{-\delta (X(t)-Y(t))} \Bigl\{ -e^{\delta (X(t)-Y(t))}\left[e^{-\delta c(t)} +\lambda(t)e^{-\delta \ell(t)}\right]  \\ \nonumber
&& +\frac{1}{\delta}(\varrho(t)+\lambda(t))+ r_r(t)X(t)+\langle\theta(t),\hat{\mu}(t)\rangle-c(t)-p(t) +h(t,z_1,z_2,z_3,\upsilon(\cdot)) \\ \nonumber
   && -\frac{\delta^2}{2} [(z_1-\theta_1b_r(t,T))^2+(z_2-(\theta_1(t)+\theta_2(t))\sigma_I(t))^2 +(z_3-\theta_3(t)\sigma_S(t))^2] \\ \nonumber
&& + \frac{1}{\delta}\int_{\mathbb{R}}[1-\delta(\langle\theta(t),\hat{\gamma}(t,T,z)\rangle-\upsilon(t,z)) -\exp(-\delta(\langle\theta(t),\hat{\gamma}(t,T,z)\rangle-\upsilon(t,z)))]\nu(dz)\Bigl\}\,.
\end{eqnarray}
Therefore, the process $\mathcal{R}_1$ is a local super-martingale if the drift process $A(t)$ is non-positive. This holds true if the generator $h$ is defined as follows

\begin{eqnarray}\label{generatorinf}
 && h(t,z_1,z_2,z_3,\upsilon) \\ \nonumber
&=& \inf_{c}\left\{e^{\delta (X(t)-Y(t))}\cdot e^{-\delta c(t)}+c(t)\right\} + \inf_{p}\left\{\lambda(t)e^{\delta (X(t)-Y(t))}\cdot e^{-\delta \ell(t)}+p(t)\right\}  \\ \nonumber
&& -\frac{1}{\delta}(\varrho(t)+\lambda(t))- r_r(t)X(t) +\inf_{\theta}\Bigl\{\frac{\delta}{2}\Bigl[\bigl|\theta_1(t)b_r(t,T)-\bigl(z_1+\frac{\varphi_1(t)}{\delta}\bigl)\bigl|^2 \\ \nonumber
   &&  +\bigl|(\theta_1(t)+\theta_2(t))\sigma_I(t)-\bigl(z_2+\frac{\varphi_2(t)}{\delta}\bigl)\bigl|^2 + \bigl|\theta_3(t)\sigma_S(t)-\bigl(z_3+\frac{\varphi_3(t)}{\delta}\bigl)\bigl|^2\Bigl] \\ \nonumber
&& +\frac{1}{\delta}\int_{\mathbb{R}}[\exp(\delta(\upsilon(t,z)-\langle\theta(t),\hat{\gamma}(t,T,z)\rangle))-1- \delta(\upsilon(t,z)-\langle\theta(t),\hat{\gamma}(t,T,z)\rangle)]\nu(dz)\Bigl\} \\ \nonumber
&& -(\varphi_1(t)z_1+\varphi_2(t)z_2+\varphi_3(t)z_3) -\frac{1}{2\delta}\left(\varphi_1^2(t) +\varphi_2^2(t)+\varphi_3^2(t)\right)\,,
\end{eqnarray}
provided that $\frac{1}{\delta}\int_{\mathbb{R}}[\exp(\delta(\upsilon(t,z)-\langle\theta(t),\hat{\gamma}\rangle))-1- \delta(\upsilon(t,z)-\langle\theta(t),\hat{\gamma}\rangle)]\nu(dz)$ is finite, for any $\theta\in\mathcal{A}$. Due to the boundedness of the associated parameters, for any $z_1,z_2,z_3\in\mathbb{R}$, the generator $h(t,z_1,z_2,z_3,\upsilon)$ is almost surely finite. Solving the three minimization problems in \eqref{generatorinf}, leads to
$$
c^*(t)=X^{(\theta^*,c^*,p^*)}(t)-Y(t)+\frac{1}{\delta}\ln\delta; \ \ \ \ p^*(t)=\eta(t)\Bigl[\frac{1}{\delta}\ln\Bigl(\frac{\delta\lambda(t)}{\eta(t)}\Bigl)-Y(t)\Bigl]
$$
and $\theta^*(t)$ in \eqref{optimalstrategy}, where $X^{(\theta^*,c^*,p^*)}$ is the wealth process associated to $(\theta^*,c^*,p^*)$ and $Y$ is part of the solution $(Y,Z_1,Z_2,Z_3, \Upsilon)$ of the BSDE with jumps \eqref{bsde}, with terminal condition $Y(T)=0$ and the generator $h$ given by \eqref{generatorh}.

To prove the super-martingale property of $\mathcal{R}_1^{(\theta,c,p)}$, we consider a function $\Psi(t)=e^{-\delta X(t)}$. Applying the generalized It\^o's formula and the dynamics of $X(t)$ in \eqref{wealth}, we have

\begin{eqnarray*}
  d\Psi(t) &=& \Psi(t)\Bigl\{\bigl[-\delta[r_r(t)X(t)+\langle\theta(t), \hat{\mu}(t)\rangle-c(t)-p(t)] \\
   && +\frac{\delta}{2}\left[\theta_1^2(t)b_r^2(t,T)+(\theta_1(t)+\theta_2(t))^2 \sigma_I^2(t) +\theta_3^2(t)\sigma_S^2(t)\right] \\
   &&  +\int_{\mathbb{R}}[\exp(-\delta \langle\theta(t), \hat{\gamma}(t,T,z)\rangle)-1+\delta\langle\theta(t), \hat{\gamma}(t,T,z)\rangle]\nu(dz)\bigl]dt \\
&& -\delta\theta_1(t) b_r(t,T)dW_r(t)-\delta (\theta_1(t)+\theta_2(t))\sigma_I(t)dW_I(t) -\delta\theta_3(t)\sigma_S(t)dW_s(t) \\
&& + \int_{\mathbb{R}}[\exp(-\delta \langle\theta(t), \hat{\gamma}(t,T,z)\rangle)-1]\tilde{N}(dt,dz)\Bigl\}\,.
\end{eqnarray*}
Therefore,

\begin{eqnarray}\nonumber
  \Psi(t) &=& \Psi(0)\mathcal{E}\Bigl( -\int_0^t\delta\theta_1(t) b_r(t,T)dW_r(t)-\int_0^t\delta (\theta_1(t)+\theta_2(t))\sigma_I(t)dW_I(t)  \\ \label{stochastice}
   && -\int_0^t\delta\theta_3(t)\sigma_S(t)dW_s(t) + \int_0^t\int_{\mathbb{R}}[\exp(-\delta \langle\theta(t), \hat{\gamma}(t,T,z)\rangle)-1]\tilde{N}(dt,dz)\Bigl)e^{K(t)}\,,
\end{eqnarray}
where $\mathcal{E}(M)$ denotes the stochastic exponential of $M$ and
\begin{eqnarray*}
  K(t) &=& \int_0^t\bigl[-\delta[r_r(s)X(s)+\langle\theta(s), \hat{\mu}(s)\rangle-c(s)-p(s)] \\
   && +\frac{\delta}{2}\left[\theta_1^2(s)b_r^2(s,t)+(\theta_1(s)+\theta_2(s))^2 \sigma_I^2(s) +\theta_3^2(s)\sigma_S^2(s)\right] \\
   &&  +\int_{\mathbb{R}}[\exp(-\delta \langle\theta(s), \hat{\gamma}(s,t,z)\rangle)-1+\delta\langle\theta(s), \hat{\gamma}(s,t,z)\rangle]\nu(dz)\bigl]ds\,.
\end{eqnarray*}
Hence, $K(t)$ is a bounded process due to the boundedness of the associated parameters and that the strategy $(c(t),p(t),\theta(t))\in\mathcal{C}\times\mathcal{D}\times\mathcal{Q}$. Furthermore, thanks to the boundedness of the associated parameters and $\exp(-\delta \langle\theta(t), \hat{\gamma}(t,T,z)\rangle)-1>-1$, the local martingale process
\begin{eqnarray*}
M(t) &:=& -\int_0^t\delta\theta_1(t) b_r(t,T)dW_r(t)-\int_0^t\delta (\theta_1(t)+\theta_2(t))\sigma_I(t)dW_I(t) \\
&& -\int_0^t\delta\theta_3(t)\sigma_S(t)dW_s(t) +\int_0^t\int_{\mathbb{R}}[\exp(-\delta \langle\theta(t), \hat{\gamma}(t,T,z)\rangle)-1]\tilde{N}(dt,dz)
\end{eqnarray*}
satisfy the $BMO$-martingale property. Then,  by {\it Kazamaki's criterion} (Morlais \cite{morlais2009}, Lemma 2.), the stochastic exponential $\mathcal{E}$ in \eqref{stochastice}, is a true martingale. Hence $\Psi$ is uniformly integrable. Then we can conclude that $\mathcal{R}_1^{(\theta,c,p)}$ is a super-martingale, i.e.,
$$
\mathcal{R}_1^{(\theta,c,p)}(0)\geq\mathbb{E}[\mathcal{R}_1^{(\theta,c,p)}(T)]\,.
$$
On the other hand, $A(t)\equiv0$, for the strategy $(\theta^*,c^*,p^*)$, hence $\mathcal{R}_1^{(\theta^*,c^*,p^*)}$ is a true martingale. Therefore, \eqref{optimalwealth} hold, which completes the proof.
\begin{flushright}
$\square$
\end{flushright}
~~\\

\subsection{The power utility case}
~~~\\

Let $\pi(t):=(\pi_1(t),\pi_2(t),\pi_3(t))$, $t\in[0,T]$ be the vector of the portfolio weights invested in $P_r^*(t,T)$, $B_r^*(t)$ and $S(t)$ respectively. Define the relative consumption rate $\xi(t)$ and the relative premium insurance rate $\zeta(t)$ by their fraction of the total wealth, i.e., $\xi(t):=\frac{c(t)}{X(t)}$ and $\zeta(t):=\frac{p(t)}{X(t)}$. We suppose that the strategy $(\pi(t),\xi(t),\zeta(t))$ satisfies the integrability condition similar to \eqref{integrabilitycp} for $(\theta(t), c(t, p(t)))$. Define $\mathcal{A}^\S$ as the admissible strategy for $(\pi(t),\xi(t),\zeta(t))$. Then, the wealth process $X(t)$ becomes

\begin{eqnarray}\label{wealthp}
  dX(t) &=& X(t)\Bigl\{[r_r(t)+\langle\pi(t),\hat{\mu}(t)\rangle-\xi(t)-\zeta(t)]dt+\pi_1(t)b_r(t,T)dW_r(t) \\ \nonumber
   && +(\pi_1(t)+\pi_2(t))\sigma_I(t)dW_I(t) +\pi_3(t)\sigma_S(t)dW_S(t) \\ \nonumber
   &&  +\int_{\mathbb{R}}\langle\pi(t),\hat{\gamma}(t,T,z)\rangle\tilde{N}(dt,dz)\Bigl\}\,,
\end{eqnarray}
which gives the following solution

\begin{eqnarray*}
  X(T) &=& x\exp\Bigl\{\int_0^T\Bigl[r_r(t)+\langle\pi(t),\hat{\mu}(t)\rangle-\xi(t)-\zeta(t) -\frac{1}{2}[\pi_1^2(t)b_r^2(t,T) \\
   && +(\pi_1(t)+\pi_2(t))^2\sigma_I^2(t)+\pi_3^2(t)\sigma_S^2(t)]  \\
&& +\int_{\mathbb{R}}[\ln(1+\langle\pi(t),\hat{\gamma}(t,T,z)\rangle)-\langle\pi(t),\hat{\gamma}(t,T,z)\rangle]\nu(dz)\Bigl]dt \\
&& +\int_0^T\pi_1(t)b_r(t,T)dW_r(t)  +\int_0^T(\pi_1(t)+\pi_2(t))\sigma_I(t)dW_I(t)  \\
&& +\int_0^T\pi_3(t)\sigma_S(t)dW_S(t) +\int_0^T\int_{\mathbb{R}}\ln(1+\langle\pi(t),\hat{\gamma}(t,T,z)\rangle)\tilde{N}(dt,dz)\Bigl\}\,.
\end{eqnarray*}

Consider the following utility function
\begin{equation}\label{powerutility}
U(x)=\frac{x^{\kappa}}{\kappa}\,, \ \ \ \ \kappa\in(-\infty,1)\setminus\{0\}\,.
\end{equation}
The functional \eqref{functional0} can be written as
\begin{eqnarray}\nonumber
  \mathcal{J}(t) &=&  \mathbb{E}_{t,x}\Bigl[\frac{1}{\kappa}\int_t^T
e^{-\int_t^s(\varrho(u)+\lambda(u))du}\left[(\xi(s))^{\kappa} +\lambda(s)\left(1+\frac{\zeta(s)}{\eta(s)}\right)^{\kappa}\right](X(s))^{\kappa}ds \\ \label{functionalpow}
   && \ \ \ \ \ \  +e^{-\int_t^T(\varrho(u)+\lambda(u))du}\frac{(X(T))^{\kappa}}{\kappa}\mid\mathcal{F}_t\Bigl]\,.
\end{eqnarray}

Define a function $B(t)$ as
\begin{eqnarray*}
  B(t) &=& \int_0^t\Bigl[r_r(s)+\langle\pi(s),\hat{\mu}(s)\rangle-\xi(s)-\zeta(s) -\frac{1}{2}[\pi_1^2(s)b_r^2(s,t) \\
   && +(\pi_1(s)+\pi_2(s))^2\sigma_I^2(s)+\pi_3^2(s)\sigma_S^2(s)]  \\
&& +\int_{\mathbb{R}}[\ln(1+\langle\pi(s),\hat{\gamma}(s,t,z)\rangle)-\langle\pi(s),\hat{\gamma}(s,t,z)\rangle]\nu(dz)\Bigl]ds \\
&& +\int_0^t\pi_1(s)b_r(s,t)dW_r(s)  +\int_0^t(\pi_1(s)+\pi_2(s))\sigma_I(s)dW_I(s)  \\
&& +\int_0^t\pi_3(s)\sigma_S(s)dW_S(s) +\int_0^t\int_{\mathbb{R}}\ln(1+\langle\pi(s),\hat{\gamma}(s,t,z)\rangle)\tilde{N}(ds,dz)\,.
\end{eqnarray*}
Then the wealth process can be written  as $X(t)=xe^{B(t)}$.

The main result of this subsection is given by the following theorem.

\begin{theorem}\label{theorem-exponential}
Suppose that the utility function is given by \eqref{powerutility}. Then, the optimal value function is given by
\begin{equation}\label{optimalpower}
V(t,x)=\frac{x^\kappa}{\kappa}e^{Y(t)}\,,
\end{equation}
where $Y$ is part of the solution $(Y,Z_1,Z_2,Z_3,\Upsilon)$ of the following BSDE with jumps
\begin{eqnarray}\label{bsdepower}
  dY(t) &=& -h_1(t, Y(t), Z_1(t), Z_2(t), Z_3(t), \Upsilon(t,\cdot))dt+Z_1(t)dW_r(t) \\ \nonumber
   && +Z_2(t)dW_I(t)+Z_3(t)dW_S(t)+\int_{\mathbb{R}}\Upsilon(t,z)\tilde{N}(dt,dz); \\ \nonumber
  Y(T) &=& 0
\end{eqnarray}
with the generator

\begin{eqnarray}\label{generatorpower}
 && h_1(t, y, z_1, z_2, z_3, \upsilon(\cdot) \\ \nonumber
&=& \left\{\frac{1}{\kappa}\Bigl(1+ \lambda(t)\Bigl(\frac{\eta(t)}{\lambda(t)}\Bigl)^{-\frac{\kappa}{1-\kappa}}\Bigl) -\Bigl(1+\eta(t) \Bigl(\frac{\eta(t)}{\lambda(t)}\Bigl)^{-\frac{1}{1-\kappa}}\Bigl) \right\} e^{-\frac{1}{1-\kappa}y} \\ \nonumber
   && -\frac{1}{\kappa}(\varrho(t)+\lambda(t))+r_r(t) +\inf_{\pi}\Bigl\{\frac{\kappa-1}{2}\Bigl[\bigl|\pi_1(t)b_r(t,T)+\frac{z_1+\varphi_1(t)}{\kappa-1}\bigl|^2 \\ \nonumber
   &&  +\bigl|(\pi_1(t)+\pi_2(t))\sigma_I(t)+\frac{z_2+\varphi_2(t)}{\kappa-1}\bigl|^2 + \bigl|\pi_3(t)\sigma_S(t)+\frac{z_3+\varphi_3(t)}{\kappa-1}\bigl|^2\Bigl] \\ \nonumber
&&\,+ \int_{\mathbb{R}}[(1+\langle\pi(t),\hat{\gamma}(t,T,z)\rangle)^{\kappa} e^{\kappa\upsilon(t,z)}-1- \kappa\langle\pi(t),\hat{\gamma}(t,T,z)\rangle-\upsilon(t,z)]\nu(dz)\Bigl\} \\ \nonumber
   && -\frac{1}{2(\kappa-1)}\bigl[(z_1+\varphi_1(t))^2 +(z_2+\varphi_2(t))^2 +(z_3+\varphi_3(t))^2\bigl] +\frac{1}{\kappa}\bigl(z_1^2+z_2^2+z_3^3\bigl)\,.
\end{eqnarray}
Moreover, the optimal strategy $(\pi^*(t),\xi^*(t),\zeta^*(t))$ is given by
$$
\xi^*(t)=e^{-\frac{1}{1-\kappa}Y(t)}, \ \ \ \ \zeta^*(t)=\eta(t)\left[\Bigl(\frac{\eta(t)}{\lambda(t)}\Bigl)^{-\frac{1}{1-\kappa}} e^{-\frac{1}{1-\kappa}Y(t)}-1\right]
$$
and
\begin{eqnarray}\nonumber
  \pi^*(t) &=& \inf_{\pi}\Bigl\{\frac{\kappa-1}{2}\Bigl[\bigl|\pi_1(t)b_r(t,T)+\frac{z_1+\varphi_1(t)}{\kappa-1}\bigl|^2 \\ \label{optimalstrategypower}
   &&  +\bigl|(\pi_1(t)+\pi_2(t))\sigma_I(t)+\frac{z_2+\varphi_2(t)}{\kappa-1}\bigl|^2 + \bigl|\pi_3(t)\sigma_S(t)+\frac{z_3+\varphi_3(t)}{\kappa-1}\bigl|^2\Bigl] \\ \nonumber
&&\,+ \int_{\mathbb{R}}[(1+\langle\pi(t),\hat{\gamma}(t,T,z)\rangle)^{\kappa} e^{\kappa\upsilon(t,z)}-1- \kappa\langle\pi(t),\hat{\gamma}(t,T,z)\rangle-\upsilon(t,z)]\nu(dz)\Bigl\}\,.
\end{eqnarray}

\end{theorem}

Note that the generator $h_1$ in \eqref{generatorpower} has an exponential growth in $Y$. However, due to the boundedness of the associated parameters, it satisfies the monotonicity condition, i.e., there exists a constant $K\geq0$ such that  $y(h_1(t, y, z_1, z_2, z_3, \upsilon(\cdot)-h_1(t, 0, z_1, z_2, z_3, \upsilon(\cdot))\leq K|y|^2$. Moreover, it can be seen that the conditions \eqref{boundedness}-\eqref{growthcondition} are satisfied. Then, by (Briand and Hu \cite{briand} and Morlais \cite{morlais2009}), the BSDE with jumps \eqref{bsdepower} has a unique solution $(Y,Z_1,Z_2,Z_3,\Upsilon)\in\mathbb{S}(\mathbb{R})\times\mathbb{H}^2(\mathbb{R})\times\mathbb{H}^2(\mathbb{R}) \times\mathbb{H}^2(\mathbb{R})\times\mathbb{H}^2_\nu(\mathbb{R})$.\\

Similar to Corollary \ref{corr-power},  the optimal investment strategy $\pi^*(t)=(\pi^*_1(t),\pi^*_2(t),\pi^*_3(t))$ for the special case of not having jumps $(\nu=0)$,  is given by the following corollary.

\begin{corollary}
Assume that $\nu=0$, then the optimal portfolio strategy $(\pi^*_1(t),\pi^*_2(t),\pi^*_3(t))$, for all $t\in[0,T]$ is given by
\begin{eqnarray*}
  \pi_1^*(t) &=& \frac{1}{1-\kappa}\Bigl[\frac{\tilde{A}(t,T)-r_r(t)-\mu_I(t)}{ b_r^2(t,T)}+\frac{Z_1(t)}{b_r(t,T)}\Bigl] \\
  \pi_2^*(t) &=& \frac{1}{1-\kappa}\Bigl[\Bigl(\frac{1}{\sigma_I^2(t)}+\frac{1}{b_r^2(t,T)}\Bigl)\mu_I(t) -\frac{\tilde{A}(t,T)-r_r(t)}{b_r^2(t,T)} +\frac{Z_2(t)}{\sigma_I(t)} -\frac{Z_1(t)}{b_r(t,T)}\Bigl] \\
  \pi_3^*(t) &=& \frac{1}{1-\kappa}\Bigl[\frac{\mu_S(t)-r_r(t)}{\sigma_S^2(t)}+\frac{Z_3(t)}{\sigma_S(t)}\Bigl]\,,
\end{eqnarray*}
where $(Z_1(t), Z_2(t), Z_3(t))$ is part of the solution $(Y, Z_1, Z_2, Z_3)$ of the following BSDE.
\begin{eqnarray*}
  dY(t) &=& -h_1(t, Y(t), Z_1(t), Z_2(t), Z_3(t), \pi^*(t), \xi^*(t),\zeta^*(t))dt+Z_1(t)dW_r(t) \\ \nonumber
   && +Z_2(t)dW_I(t)+Z_3(t)dW_S(t); \\ \nonumber
  Y(T) &=& 0\,.
\end{eqnarray*}
\end{corollary}
~~\\

\noindent {\it Proof of Theorem \ref{theorem-exponential}.}

Consider the process
\begin{eqnarray}\nonumber
  \mathcal{R}_2(t) &=&  \frac{1}{\kappa}\int_0^t
e^{-\int_0^s(\varrho(u)+\lambda(u))du}\Bigl[(\xi(s))^{\kappa} +\lambda(s)\Bigl(1+\frac{\zeta(s)}{\eta(s)}\Bigl)^{\kappa}\Bigl](X(s))^{\kappa}ds \\ \nonumber
   && \ \ \ \ \ \  +e^{-\int_0^t(\varrho(u)+\lambda(u))du}\frac{(X(t))^{\kappa}}{\kappa}e^{Y(t)}\,,
\end{eqnarray}
with initial condition $\mathcal{R}_2(0)= \frac{x^{\kappa}}{\kappa}e^{Y(0)}$\,. Applying the generalized It\^o's formula, we obtain

\begin{eqnarray}\nonumber
  d\mathcal{R}_2(t) &=& e^{-\int_0^t(\varrho(u)+\lambda(u))du}(X(t))^{\kappa}e^{Y(t)}\Bigl\{\Bigl[\frac{1}{\kappa}e^{-Y(t)} \Bigl((\xi(t))^{\kappa} +\lambda(t)\Bigl(1+\frac{\zeta(t)}{\eta(t)}\Bigl)^{\kappa}\Bigl) \\ \nonumber
   && -\frac{1}{\kappa}(\varrho(t)+\lambda(t))+r_r(t)+\langle\pi(t),\hat{\mu}(t)\rangle -\xi(t)-\zeta(t)-h_1(t,y,z_1,z_2,z_3,\upsilon) \\ \nonumber
   && +\frac{1}{2}(\kappa-1)[\pi_1^2(t)b_r^2(t,T)+(\pi_1(t)+\pi_2(t))^2\sigma_I^2(t) +\pi_3^2(t)\sigma_S^2(t)] \\ \nonumber
&& +\frac{1}{2\kappa}(z_1^2+z_2^2+z_3^2)+ \pi_1(t)b_r(t,T)z_1+(\pi_1(t)+\pi_2(t))\sigma_I(t)z_2 +\pi_3(t)\sigma_S(t)z_3 \\ \nonumber
&& +\frac{1}{\kappa}\int_{\mathbb{R}}[(1+\langle\pi(t),\hat{\gamma}(t,T,z)\rangle)^{\kappa} e^{\kappa\upsilon(t,z)}-1- \kappa\langle\pi(t),\hat{\gamma}(t,T,z)\rangle-\upsilon(t,z)]\nu(dz)\Bigl]dt \\ \nonumber
&& +\frac{1}{\kappa}[(z_1+\kappa\pi_1(t)b_r(t,T))dW_r(t) +(z_2+\kappa(\pi_1(t)+\pi_2(t))\sigma_I(t))dW_I(t) \\ \nonumber
&& +(z_3+\kappa\pi_3(t)\sigma_S(t))dW_S(t)] + \frac{1}{\kappa}\int_{\mathbb{R}}[(1+\langle\pi(t),\hat{\gamma}(t,T,z)\rangle)^{\kappa} e^{\kappa\upsilon(t,z)}-1  \\ \label{r2process}
&& +\kappa\ln(1+\langle\pi(t),\hat{\gamma}(t,T,z)\rangle) -\upsilon(t,z)]\tilde{N}(dz,dt)\Bigl\}\,.
\end{eqnarray}
Note that similar to the exponential case, we can easily see that the process $\mathcal{R}_2$ is a local super-martingale if the generator $h_1$ is given by

\begin{eqnarray*}
&& h_1(t,y,z_1,z_2,z_3,\upsilon) \\
&=& \inf_{\xi}\Bigl\{\frac{1}{\kappa}e^{-Y(t)} (\xi(t))^{\kappa}-\xi(t)\Bigl\} +\inf_{\zeta}\Bigl\{ \frac{1}{\kappa}\lambda(t)e^{-Y(t)} \Bigl(1+\frac{\zeta(t)}{\eta(t)}\Bigl)^{\kappa}-\zeta(t)\Bigl\}  \\
   && + \inf_{\pi}\Bigl\{\frac{\kappa-1}{2}\Bigl[\bigl|\pi_1(t)b_r(t,T)+\frac{z_1+\varphi_1(t)}{\kappa-1}\bigl|^2 \\ \nonumber
   &&  +\bigl|(\pi_1(t)+\pi_2(t))\sigma_I(t)+\frac{z_2+\varphi_2(t)}{\kappa-1}\bigl|^2 + \bigl|\pi_3(t)\sigma_S(t)+\frac{z_3+\varphi_3(t)}{\kappa-1}\bigl|^2\Bigl] \\ \nonumber
&&\,+ \int_{\mathbb{R}}[(1+\langle\pi(t),\hat{\gamma}(t,T,z)\rangle)^{\kappa} e^{\kappa\upsilon(t,z)}-1- \kappa\langle\pi(t),\hat{\gamma}(t,T,z)\rangle-\upsilon(t,z)]\nu(dz)\Bigl\}  \\
   &&  -\frac{1}{2(\kappa-1)}\bigl[(z_1+\varphi_1(t))^2 +(z_2+\varphi_2(t))^2 +(z_3+\varphi_3(t))^2\bigl] \\
&& +\frac{1}{\kappa}\bigl[z_1^2+z_2^2+z_3^3\bigl] -\frac{1}{\kappa}(\varrho(t)+\lambda(t))+r_r(t)\,.
\end{eqnarray*}
Solving the three minimization problems, provided that the associated parameters are bounded $\mathcal{F}_t$-predictable, we obtain the candidate optimal strategy
$$
\xi^*(t)=e^{-\frac{1}{1-\kappa}Y(t)}, \ \ \ \ \zeta^*(t)=\eta(t)\left[\Bigl(\frac{\eta(t)}{\lambda(t)}\Bigl)^{-\frac{1}{1-\kappa}} e^{-\frac{1}{1-\kappa}Y(t)}-1\right]
$$
and $\pi^*(t)$ in \eqref{optimalstrategypower}. Where  $Y$ is part of the solution $(Y,Z_1,Z_2,Z_3, \Upsilon)$ of the BSDE with jumps \eqref{bsdepower}, with terminal condition $Y(T)=0$ and the generator $h_1$ given by \eqref{generatorpower}.

To prove the super-martingale property, we consider the following function $(X(t))^{\kappa}$. Applying the generalized It\^o's formula, we have
\begin{eqnarray*}
  d(X(t))^{\kappa} &=& (X(t))^{\kappa}\Bigl\{\Bigl[r_r(t)+\langle\pi(t),\hat{\mu}(t)\rangle-\xi(t)-\zeta(t) +\frac{\kappa}{2}(\kappa-1)[\pi_1^2(t)b_r^2(t,T) \\
   && +(\pi_1(t)+\pi_2(t))^2\sigma_I^2(t)+\pi_3^2(t)\sigma_S^2(t)]  \\
&& +\int_{\mathbb{R}}[1-\kappa\ln(1+\langle\pi(t),\hat{\gamma}(t,T,z)\rangle)- (1+\langle\pi(t),\hat{\gamma}(t,T,z)\rangle)^{-\kappa}]\nu(dz)\Bigl]dt \\
   &&  +\kappa\pi_1(t)b_r(t,T)dW_r(t)  +\kappa(\pi_1(t)+\pi_2(t))\sigma_I(t)dW_I(t)  \\
&& +\kappa\pi_3(t)\sigma_S(t)dW_S(t) +\kappa\int_{\mathbb{R}}\ln(1+\langle\pi(t),\hat{\gamma}(t,T,z)\rangle)\tilde{N}(dt,dz)\Bigl\}\,.
\end{eqnarray*}
Hence
\begin{eqnarray}\nonumber
 (X(t))^{\kappa}  &=& (X(0))^{\kappa}\mathcal{E}\Bigl(\kappa\pi_1(t)b_r(t,T)dW_r(t)  +\kappa(\pi_1(t)+\pi_2(t))\sigma_I(t)dW_I(t)  \\
&& +\kappa\pi_3(t)\sigma_S(t)dW_S(t) +\kappa\int_{\mathbb{R}}\ln(1+\langle\pi(t),\hat{\gamma}(t,T,z)\rangle)\tilde{N}(dt,dz)\Bigl)e^{Q(t)}\,,
\end{eqnarray}
where
\begin{eqnarray*}
  Q(t) &=& \int_0^t\Bigl[r_r(s)+\langle\pi(s),\hat{\mu}(s)\rangle-\xi(s)-\zeta(s) +\frac{\kappa}{2}(\kappa-1)[\pi_1^2(s)b_r^2(s,t) \\
   && +(\pi_1(s)+\pi_2(s))^2\sigma_I^2(s)+\pi_3^2(s)\sigma_S^2(s)]  \\
&& +\int_{\mathbb{R}}[1-\kappa\ln(1+\langle\pi(s),\hat{\gamma}(s,t,z)\rangle)- (1+\langle\pi(s),\hat{\gamma}(s,t,z)\rangle)^{-\kappa}]\nu(dz)\Bigl]ds\,.
\end{eqnarray*}
Hence, $Q(t)$ is a bounded process due to the boundedness of the associated parameters and the fact that the strategy $(\xi(t),\zeta(t),\pi(t))\in\mathcal{C}\times\mathcal{D}\times\mathcal{Q}$. Moreover, using similar arguments of a $BMO$-martingale property as in the proof of Theorem \ref{theorem-power}, we can easily see that $(X(t))^{\kappa}$ is uniformly integrable. Then $\mathcal{R}_2^{(\theta,c,p)}$ is a super-martingale, i.e.,
$$
\mathcal{R}_2^{(\theta,c,p)}(0)\geq\mathbb{E}[\mathcal{R}_2^{(\theta,c,p)}(T)]\,.
$$
On the other hand, the drift process in \eqref{r2process} is equal to zero for the strategy $(\theta^*,c^*,p^*)$, hence $\mathcal{R}_1^{(\theta^*,c^*,p^*)}$ is a true martingale. Therefore, \eqref{optimalpower} hold.
\begin{flushright}
$\square$
\end{flushright}
\begin{remark}
We point out that, when there is no inflation and jumps in the model, the results obtained in Theorems \ref{theorem-power} and \ref{theorem-exponential} relate on the results in Cheridito and Hu \cite{cheridito}. Similar results have been obtained by Xing \cite{xing} for the Espein-Zin utility type. Moreover, when there is no consumption and life insurance rates, these results are similar to those in Hu {\it et. al.} \cite{hu} for the diffusion case and Morlais \cite{morlais2009} in the jump-diffusion case.
\end{remark}

\section{Conclusion}

In this paper, we solved an optimal investment, consumption and life insurance problem using the BSDE techniques. We considered the presence of inflation-linked asset, which normally help the investors to manage the inflation risks that in general are not completely observable. Under jump diffusion market, we derived the optimal strategy for the exponential and power utility functions. This work extends, for instance, the paper by Cheridito and HU \cite{cheridito}, by allowing the presence of inflation risks, life insurance and jumps in the related assets. Furthermore, it appears as an alternative approach to the dynamic programming approach applied in Han and Hung \cite{Han}, were a similar problem was considered under a stochastic differential utility. We noted that the generator of the associated BSDE is of quadratic growth in the controls $z_1,\,z_2,\, z_3$ and exponential in $\upsilon(\cdot)$,  the similar BSDEs with jumps that the existence and uniqueness results have been proved by Morlais \cite{morlais2009}, \cite{morlais2010}. Furthermore, we derived the explicit solutions for the optimal portfolio for a special case without jumps.

\subsection*{Acknowledgment}

We would like to express our deep gratitude to the NRF Project No: CSUR 90313, the University of Pretoria and the MCTESTP Mozambique for their support. We also wish to thank the four anonymous referees for their comments and critics during the reviewing process

\end{document}